\documentclass[11pt]{article}
\usepackage{latexsym,a4wide}
\usepackage{amssymb}
\usepackage{amsmath}
\usepackage{color}
\usepackage{graphicx}
\usepackage{amsthm}
\usepackage{indentfirst}
\allowdisplaybreaks

\newcommand \bel {\begin{equation}\label}
\newcommand \ee {\end{equation}}
\newcommand \be {\begin{equation}}

\newcommand \RR {\mathbb R}

\newcommand \del \partial

\newcommand \bei {\begin{itemize}}
\newcommand \eei {\end{itemize}}

\newtheorem{theorem}{\color{black}\indent Theorem}[section]

\newtheorem{proposition}{\color{black}\indent Proposition}[section]

\newtheorem{remark}{\color{black}\indent Remark}[section]

\begin{document}
\large
\title{On the explicit blowup solutions for $3$D incompressible Magnetohydrodynamics equations}
\author{
{\sc Weiping Yan}\thanks{School of Mathematics, Xiamen University, Xiamen 361000, P.R. China. Email: yanwp@xmu.edu.cn.}
\thanks{Laboratoire Jacques-Louis Lions, Sorbonne Université, 4, Place Jussieu, 75252 Paris, France.}
}
\date{July 3, 2018}

\maketitle

\begin{abstract}  
This paper concerns with the explicit blowup phenomenon for $3$D incompressible MHD equations in $\RR^3$.
More precisely, we find two family of explicit blowup solutions for $3$D incompressible MHD equations in $\RR^3$.
One family of solutions admit the smooth initial data, and the initial data of another family of solutions are not smooth. The energy of those solutions is infinite.
Moreover, our results tell us that the blowup phenomenon of $3$D incompressible MHD can only take place in the velocity field of the fluid, but no blowup for the magnetic field.
\end{abstract}



\section{Introduction and main results} 
\setcounter{equation}{0}

The $3$D incompressible Magnetohydrodynamics equations (MHD) describes the dynamics of electrically conducting fluids arising from plasmas or some other physical phenomena. 
We are concerned with the blowup phenomena of smooth solutions to the MHD equations in $\RR^3$:
\bel{E1-1}
\aligned
&\del_t\textbf{v}+\textbf{v}\cdot\nabla \textbf{v}+\nabla P=\nu\triangle \textbf{v}+(\nabla\times\textbf{H})\times\textbf{H},\\
&\del_t\textbf{H}=\nu\triangle\textbf{H}+\nabla\times(\textbf{v}\times\textbf{B}),\\
&\nabla\cdot\textbf{v}=0,\quad \nabla\cdot\textbf{H}=0,
\endaligned
\ee
where $(t,x)\in\RR\times\RR^3$, $\textbf{v}$ denotes the $3$D velocity field of the fluid, $P$ stands for the pressure in the fluid, $\textbf{H}$ is the the magnetic field,
$\nu\geq0$ denotes the viscosity constant.
The divergence free condition in second equations of (\ref{E1-1}) guarantees the incompressibility of the fluid. 
In particularly, if $\nu=0$, equations (\ref{E1-1}) is called the ideal MHD equations.


It is easy to check that solutions of $3$D incompressible MHD equations (\ref{E1-1}) admits the scaling invariant property, that is, 
let $(\textbf{v},\textbf{H},P)$ be a solution of (\ref{E1-1}),  then for any constant $\lambda>0$ and $\alpha\in\RR$, the functions
$$
\aligned
&\textbf{v}_{\lambda,\alpha}(t,x)=\lambda^{\alpha}\textbf{v}(\lambda^{\alpha+1}t,\lambda x),\\
&\textbf{H}_{\lambda,\alpha}(t,x)=\lambda^{\alpha}\textbf{H}(\lambda^{\alpha+1}t,\lambda x),\\
&P_{\lambda,\alpha}(t,x)=\lambda^{2\alpha}P(\lambda^{\alpha+1}t,\lambda x),
\endaligned
$$
are also solutions of $3$D incompressible MHD equations (\ref{E1-1}). 
Here the initial data $(\textbf{v}_0(x),\textbf{H}_0(x))$ is changed into $(\lambda^{\alpha}\textbf{v}_0(\lambda x),\lambda^{\alpha}\textbf{H}_0(\lambda x))$.

The MHD equations (\ref{E1-1}) is a combination of the Navier-Stokes equations of fluid dynamics and Maxwell's equations of electromagnetism.
As we known, the question of finite time singularity/global regularity for $3$D incomprsssible Navier-Stokes equations is the most important open problems in mathematical fluid mechanics
\cite{F}. In 1934, Leray \cite{L} showed that the $3$D incompressible Navier-Stokes equations admit global-forward-in-time weak solution of the initial value problem. After that, there are many papers concerns with the well-posedness of weak solutions or blowup of solution for this problem. One can see \cite{CKN,CP,HL,JS,NRS} for more details. 
So it is also natural important problem for the $3$D incompressible MHD equations.
There are some numerical results to approach the singularity of this kind problem \cite{GO}. The Beale-Kato-Majda's blowup criterion for incompressible MHD was obtained in \cite{CKS,CM}.
Chae \cite{Ch} excluded the scenario of the apparition of finite time singularity in the form of self-similar singularities.
In this paper, we will give two family of explicit blowup solutions for $3$D incompressible MHD equations (\ref{E1-1}). One can see that the ideal MHD equations also has the same explicit blowup solutions with equations (\ref{E1-1}). This means the viscosity constant $\nu$ can not effect the blowup phenomenon.

For the existence of explicit solutions of incompressible fluids,
Landau \cite{Lan} gave a family of explicit stationary solutions for $3$D incompressible stationary Navier-Stokes equations, those solutions now are called Landau's solutions \cite{Lan0}.
Sver\'{a}k \cite{Sv} proved that there is only a Landau solution for a nontrivial smooth solutions of $3$D incompressible stationary Navier-Stokes equations if the solution satisfying the scaling invariant 
$\textbf{v}(x)=\lambda\textbf{v}(\lambda x)$ for each constant $\lambda>0$.
Kapitanskiy \cite{K,K0} found the first nontrivial example of hidden symmetries connected with reduction of PDEs for the Navier-Stokes equations.
Constantin \cite{Con} found a class of smooth, mean zero initial data for which the solution of $3$D Euler equations becomes infinite in finite time, meanwhile, he gave an explicit formulas of solutions for the $3$D Euler equations by reducing this equations into a local conservative Riccati system in two-dimensional basic square. 
Recently, Yan \cite{Yan2} constructed two family of explicit blowup solutions for $3$D incompressible Navier-Stokes equations in $\RR^3$. One family of explicit blowup solutions admit smooth data and infinite energy, another family of explicit blowup solutions have non-smooth data.
Due to space limitations, here we can not list all of interesting results.
To the author's knowledgement, there is few result on the explicit meaningful solution of incompressible MHD equations (\ref{E1-1}) except the famous Alfv\'{e}n waves \cite{A}. 
One can see \cite{CL,HXY} for the result on nonlinear stability of Alfv\'{e}n waves.

Here we state the main result of this paper.

\begin{theorem}
\begin{itemize}
Let constant $T^*>0$ be maximal existence time and the viscosity constant  $\nu\geq0$.
The 3D incompressible MHD equations (\ref{E1-1}) admits two family of explicit blowup solutions as follows

\item
One family of explicit blowup solutions are
\bel{E1-7}
\aligned
&\textbf{v}(t,x)=\Big(v_1(t,x),v_2(t,x),v_3(t,x)\Big)^T,\quad (t,x)\in[0,T^*)\times\RR^3,\\
&\textbf{H}(t,x)=\Big(H_1(t,x),H_2(t,x),H_3(t,x)\Big)^T,\quad (t,x)\in[0,T^*)\times\RR^3,
\endaligned
\ee
where
$$
\aligned
&v_1(t,x):={ax_1\over T^*-t}+{kx_2\over x_1^2+x_2^2},\\
&v_2(t,x):={ax_2\over T^*-t}-{kx_1\over x_1^2+x_2^2},\\
&v_3(t,x):=-{2ax_3\over T^*-t},
\endaligned
$$
and
$$
\aligned
&H_1(t,x):=\bar{a}x_1+{2\bar{a}kx_2(T^*-t)\over (2a+1)(x_1^2+x_2^2)},\\
&H_2(t,x):=\bar{a}x_2-{2\bar{a}kx_1(T^*-t)\over (2a+1)(x_1^2+x_2^2)},\\
&H_3(t,x):=-2\bar{a}x_3,
\endaligned
$$
with the initial data
\bel{E1-8}
\aligned
&\textbf{v}(0,x)=\Big({ax_1\over T^*}+{kx_2\over x_1^2+x_2^2},~{ax_2\over T^*}-{kx_1\over x_1^2+x_2^2},~-{2ax_3\over T^*}\Big)^T,\\
&\textbf{H}(0,x)=\Big(\bar{a}x_1+{2\bar{a}kT^*x_2\over (2a+1)(x_1^2+x_2^2)},\bar{a}x_2-{2\bar{a}kx_1T^*\over (2a+1)(x_1^2+x_2^2)},-2\bar{a}x_3\Big)^T,
\endaligned
\ee
where constants $k,\bar{a}\in\RR/\{0\}$ and $a\in\RR/\{-{1\over2},0\}$.

\item
Another family of explicit blowup soutions are
\bel{E1-7R1}
\aligned
&\textbf{v}(t,x)=\Big(v_1(t,x),v_2(t,x),v_3(t,x)\Big)^T,\quad (t,x)\in[0,T^*)\times\RR^3,\\
&\textbf{H}(t,x)=\Big(H_1(t,x),H_2(t,x),H_3(t,x)\Big)^T,\quad (t,x)\in[0,T^*)\times\RR^3,
\endaligned
\ee
where
$$
\aligned
&v_1(t,x):={ax_1\over T^*-t}+kx_2(T^*-t)^{2a},\\
&v_2(t,x):={ax_2\over T^*-t}-kx_1(T^*-t)^{2a},\\
&v_3(t,x):=-{2ax_3\over T^*-t},
\endaligned
$$
and
$$
\aligned
&H_1(t,x):=\bar{a}x_1+{2\bar{a}kx_2x_3(T^*-t)^{2a+1}\over 4a+1},\\
&H_2(t,x):=\bar{a}x_2-{2\bar{a}kx_1x_3(T^*-t)^{2a+1}\over 4a+1},\\
&H_3(t,x):=-2\bar{a}x_3,
\endaligned
$$
with the smooth initial data
\bel{E1-8R1}
\aligned
&\textbf{v}(0,x)=\Big({ax_1\over T^*}+kx_2(T^*)^{2a},~{ax_2\over T^*}-kx_1(T^*)^{2a},~-{2ax_3\over T^*}\Big)^T,\\
&\textbf{H}(0,x)=\Big(\bar{a}x_1+{2\bar{a}kx_2x_3(T^*)^{2a+1}\over 4a+1},\bar{a}x_2-{2\bar{a}kx_1x_3(T^*)^{2a+1}\over 4a+1},-2\bar{a}x_3\Big)^T,
\endaligned
\ee
where constants $k,\bar{a}\in\RR/\{0\}$ and $a\in\RR/\{-{1\over4},0\}$.
\end{itemize}
\end{theorem}

\begin{remark}
If the $3$D velocity field of the fluid $\textbf{v}$ and the the magnetic field $\textbf{H}$ are given by (\ref{E1-7}), respectively, by direct computations, the pressure
\bel{XX1-1}
P(t,x)=-{1\over 2}\Big({a(a+1)(x_1^2+x_2^2)\over (T^*-t)^2}+{k^2\over x_1^2+x_2^2}\Big)+{ax_3^2(1-2a)\over (T^*-t)^2},
\ee
where constants $k\in\RR/\{0\}$ and $a\in\RR/\{-{1\over2},0\}$.

If the $3$D velocity field of the fluid $\textbf{v}$ and the the magnetic field $\textbf{H}$ are given by (\ref{E1-7R1}), respectively, by direct computations, the pressure
\bel{XX1-2}
\aligned
P(t,x)&={x_1^2+x_2^2\over 2}\Big(k^2(T^*-t)^{4a}-{a(a+1)\over (T^*-t)^2}-{8\bar{a}^2k^2x_3^2(T^*-t)^{2(2a+1)}\over (4a+1)^2}\Big)\\
&\quad+x_3^2\Big({a(1-2a)\over (T^*-t)^2}-{2\bar{a}^2k^2r^2(T^*-t)^{2(2a+1)}\over (4a+1)^2}\Big),
\endaligned
\ee
where constants $k,\bar{a}\in\RR/\{0\}$ and $a\in\RR/\{-{1\over4},0\}$.

From (\ref{XX1-1}), we can see if the blowup phenomenon only takes place in the velocity field, then the pressure $P$ depends on two parameters $a,k$.
But if the blowup phenomenon takes place in both the velocity field and the magnetic field, it follows from (\ref{XX1-2}) that the pressure $P$ depends on three parameters $a,k,\bar{a}$, where the parameter $\bar{a}$ comes from the magnetic field.
\end{remark}

\begin{remark}
For the velocity field of the fluid $\textbf{v}$, it follows from (\ref{E1-7}) and (\ref{E1-7R1}) that there is self-smilar singularity in $x_3$ direction, that is, $-{2ax_3\over T^*-t}$ for $a\in\RR/\{0\}$.
Moreover, on one hand, by (\ref{E1-7}), we find there is only blowup for the velocity field of the fluid $\textbf{v}$, the magnetic field does not blowup as $t\rightarrow (T^*)^-$.
On the other hand, by (\ref{E1-7R1}), there are not only blowup for velocity field of the fluid $\textbf{v}$, but also blowup for the magnetic field $\textbf{H}$ with constant $a<-{1\over 2}$ as $t\rightarrow (T^*)^-$.
\end{remark}

\begin{remark}

On one hand, it follows from (\ref{E1-7}) that
$$
\aligned
&\nabla v_1(t,x)=\Big({a\over T^*-t}-{2kx_1x_2\over(x_1^2+x_2^2)^2},~{k(x_1^2-x_2^2)\over (x_1^2+x_2^2)^2},~0\Big)^T,\\
&\nabla v_2(t,x)=\Big({k(x_1^2-x_2^2)\over (x_1^2+x_2^2)^2},~{a\over T^*-t}+{2kx_1x_2\over(x_1^2+x_2^2)^2},~0\Big)^T,\\
&\nabla v_3(t,x)=\Big(0,~0,~-{2a\over T^*-t}\Big)^T,
\endaligned
$$
which means that
$$
div(v_i)|_{x=x_0}=\infty,\quad as\quad t\rightarrow (T^*)^-,
$$
for a fixed point $x_0\in\RR^3$. There is no blowup in the magnetic field $\textbf{H}$.
Here one can see the initial data is not smooth from (\ref{E1-8}).

On the other hand, it follows from (\ref{E1-7R1}) that
$$
\aligned
&\nabla v_1(t,x)=\Big({a\over T^*-t},~k(T^*)^{2a},~0\Big)^T,\\
&\nabla v_2(t,x)=\Big(-k(T^*)^{2a},~{a\over T^*-t},~0\Big)^T,\\
&\nabla v_3(t,x)=\Big(0,~0,~-{2a\over T^*-t}\Big)^T,\\
&\nabla H_1(t,x)=\Big(\bar{a},{2\bar{a}kx_3(T^*-t)^{2a+1}\over 4a+1},{2\bar{a}kx_2(T^*-t)^{2a+1}\over 4a+1}\Big)^T,\\
&\nabla H_2(t,x)=\Big(-{2\bar{a}kx_3(T^*-t)^{2a+1}\over 4a+1},\bar{a},-{2\bar{a}kx_1(T^*-t)^{2a+1}\over 4a+1}\Big)^T,\\
&\nabla H_3(t,x)=\Big(0,0,-2\bar{a}\Big)^T,
\endaligned
$$
which means that
$$
div(v_i)|_{x=x_0}=\infty,\quad as\quad t\rightarrow (T^*)^-,
$$
and for $a<-{1\over 2}$, there is
$$
div(H_i)|_{x=x_0}=\infty,\quad as\quad t\rightarrow (T^*)^-,
$$
for a fixed point $x_0\in\RR^3$. Here one can see the initial data is smooth from (\ref{E1-8R1}). But the initial data goes to infinity as $x\rightarrow\infty$.

In conclusion, our results tell us that the blowup phenomenon of $3$D incompressible MHD (\ref{E1-1}) can only take place in the velocity field of the fluid $\textbf{v}$, but no blowup for the magnetic field $\textbf{H}$. It is easy to see the blowup solutions (\ref{E1-7}) and (\ref{E1-7R1}) independent of viscosity constant $\nu$, so our results also hold for $3$D ideal incompressible MHD.
Hence those two family of blowup solutions (\ref{E1-7}) and (\ref{E1-7R1}) are not finite energy solutions.
\end{remark}

\begin{remark}
When the magnetic field $\textbf{H}\equiv0$, equations (\ref{E1-1}) is reduced into $3$D incompressible Navier-Stokes equations. Then corresponding explicit blowup solutions given in (\ref{E1-7}) and (\ref{E1-7R1}) are also explicit blowup solutions for $3$D incompressible Navier-Stokes equations \cite{Yan2}.
\end{remark}

\section{Proof of Thereom 1.1}

We first recall a result on the existence of explicit blowup axisymmetric solutions for $3$D incompressible Navier-Stokes equations.

Let $\textbf{e}_r$, $\textbf{e}_{\theta}$ and $\textbf{e}_z$ be the cylindrical coordinate
system,
\bel{E1-0}
\aligned
&\textbf{e}_r=({x_1\over r},{x_2\over r},0)^T,\\
&\textbf{e}_{\theta}=({x_2\over r},-{x_1\over r},0)^T,\\
&\textbf{e}_z=(0,0,1)^T,
\endaligned
\ee
where $r=\sqrt{x_1^2+x_2^2}$ and $z=x_3$.

The following result is taken from \cite{Yan2}.

\begin{proposition} 
\begin{itemize}
Let $T^*>0$ be a constant,  the kinematic viscosity $\nu>0$, and $\textbf{e}_r,\textbf{e}_{\theta},\textbf{e}_z$ are defined in (\ref{E1-0}), $r=\sqrt{x_1^2+x_2^2}$ and $z=x_3$.
The 3D incompressible Navier-Stokes equations admits two family of explicit blowup axisymmetric solutions: 

\item One family of explicit blowup axisymmetric solutions are
\bel{E1-6}
\textbf{v}(t,x)=v^r(t,r,z)\textbf{e}_r+v^{\theta}(t,r,z)\textbf{e}_{\theta}+v^z(t,r,z)\textbf{e}_z,\quad (t,x)\in[0,T^*)\times\RR^3,
\ee
where
$$
\aligned
&v^r(t,r,z)={ar\over T^*-t},\\
&v^{\theta}(t,r,z)={k\over r},\\
&v^z(t,r,z)=-{2az\over T^*-t},
\endaligned
$$
where constants $a,k\in\RR/\{0\}$.

\item Another family of explicit blowup axisymmetric solutions are
\bel{E1-6r1}
\textbf{v}(t,x)=v^r(t,r,z)\textbf{e}_r+v^{\theta}(t,r,z)\textbf{e}_{\theta}+v^z(t,r,z)\textbf{e}_z,\quad (t,x)\in[0,T^*)\times\RR^3,
\ee
where
$$
\aligned
&v^r(t,r,z)={ar\over T^*-t},\\
&v^{\theta}(t,r,z)= kr(T^*-t)^{2a},\\
&v^z(t,r,z)=-{2az\over T^*-t},
\endaligned
$$
where constants $a,k\in\RR/\{0\}$.
\end{itemize}
\end{proposition}

We now derive the $3$D incompressible MHD equations (\ref{E1-1}) with axisymmetric velocity field in the cylindrical coordinate (e.g. see \cite{Lei}).
The $3$D velocity field $\textbf{v}(t,x)$ and magnetic field $\textbf{H}(t,x)$ are called axisymmetric if they can be written as
$$
\aligned
&\textbf{v}(t,x)=v^r(t,r,z)\textbf{e}_r+v^{\theta}(t,r,z)\textbf{e}_{\theta}+v^z(t,r,z)\textbf{e}_z,\\
&\textbf{H}(t,x)=H^r(t,r,z)\textbf{e}_r+H^{\theta}(t,r,z)\textbf{e}_{\theta}+H^z(t,r,z)\textbf{e}_z,\\
&P(t,x)=P(t,r,z),
\endaligned
$$
where $(v^r,v^{\theta},v^z)$, $(H^r,H^{\theta},H^z)$ and $P(t,r,z)$ do not depend on the $\theta$ coordinate.

Note that the Lorentz force term 
$$
(\nabla\times\textbf{H})\times\textbf{H}=\textbf{H}\cdot\nabla\textbf{H}-\nabla{|\textbf{H}|^2\over 2}.
$$
Then the $3$D MHD equations (\ref{E1-1}) with axisymmetric velocity field in the cylindrical coordinates can be reduced into a system as follows 
\bel{X1-2}
\del_tv^r+v^r\del_rv^r+v^z\del_zv^r-{1\over r}(v^{\theta})^2+\del_r\bar{P}
=\nu(\triangle-{1\over r^2})v^r+H^r\del_rH^r+H^z\del_zH^r-{1\over r}(H^{\theta})^2,
\ee
\bel{X1-2R0}
\del_tv^{\theta}+v^r\del_rv^{\theta}+v^z\del_zv^{\theta}+{1\over r}v^rv^{\theta}
=\nu(\triangle-{1\over r^2})v^{\theta}+H^r\del_rH^{\theta}+H^z\del_zH^{\theta}+{1\over r}H^{\theta}H^r,
\ee
\bel{X1-2R1}
\del_tv^z+v^r\del_rv^z+v^z\del_zv^z+\del_z\bar{P}=\nu\triangle v^z+H^r\del_rH^z+H^z\del_zH^z,
\ee
\bel{X1-2R2}
\del_tH^r+v^r\del_rH^r+v^z\del_zH^r=\nu(\triangle-{1\over r^2})H^r+H^r\del_rv^r+H^z\del_zv^r,
\ee
\bel{X1-2R3}
\del_tH^{\theta}+v^r\del_rH^{\theta}+v^z\del_zH^{\theta}+{1\over r}H^rv^{\theta}
=\nu(\triangle-{1\over r^2})H^{\theta}+H^r\del_rv^{\theta}+H^z\del_zv^{\theta}+{1\over r}v^rH^{\theta},
\ee
\bel{X1-2R4}
\del_tH^z+v^r\del_rH^z+v^z\del_zH^z=\nu\triangle H^z+H^r\del_rv^z+H^z\del_zv^z,
\ee
where the pressure is given by 
\bel{X1-3}
\bar{P}=P+{|\textbf{H}|^2\over 2}.
\ee

The incompressibility condition becomes
\bel{X1-4}
\aligned
&\del_r(rv^r)+\del_z(rv^z)=0,\\
&\del_r(rH^r)+\del_z(rH^z)=0.
\endaligned
\ee

The following result gives two family of explict self-similar blowup solutions for system (\ref{X1-2})-(\ref{X1-3}) with the  incompressibility condition (\ref{X1-4}).

\begin{proposition}
\begin{itemize}
Let $T^*>0$ be a constant and $\nu>0$. System (\ref{X1-2})-(\ref{X1-3}) with the incompressibility condition (\ref{X1-4}) admits two family of explicit blowup solutions:

\item One family of explicit blowup solutions takes the form
\bel{X1-5}
\aligned
&v^r(t,r,z)={ar\over T^*-t},\\
&v^{\theta}(t,r,z)={k\over r},\\
&v^z(t,r,z)=-{2az\over T^*-t},\\
&H^r(t,r,z)=\bar{a}r,\\
&H^{\theta}(t,r,z)={2k\bar{a}(T^*-t)\over (2a+1)r},\\
&H^z(t,r,z)=-2\bar{a}z,
\endaligned
\ee
where constants $\bar{a},k\in\RR/\{0\}$ and $a\in\RR/\{-{1\over2},0\}$.

\item Another family of explicit blowup solutions is
\bel{X1-6}
\aligned
&v^r(t,r,z)={ar\over T^*-t},\\
&v^{\theta}(t,r,z)=kr(T^*-t)^{2a},\\
&v^z(t,r,z)=-{2az\over T^*-t},\\
&H^r(t,r,z)=\bar{a}r,\\
&H^{\theta}(t,r,z)={2\bar{a}krz(T^*-t)^{2a+1}\over 4a+1},\\
&H^z(t,r,z)=-2\bar{a}z,
\endaligned
\ee
where constants $\bar{a},k\in\RR/\{0\}$ and $a\in\RR/\{-{1\over 4},0\}$.

\end{itemize}
\end{proposition}

\paragraph{Proof of Proposition 2.2.}
The idea of finding explicit blowup solutions for system (\ref{X1-2})-(\ref{X1-3}) with the incompressibility condition (\ref{X1-4}) comes from \cite{Yan1,Yan2,Yan0}. This is based on the observation on the structure of system (\ref{X1-2})-(\ref{X1-3}) and incompressibility condition (\ref{X1-4}). We notice that if the magnetic field $\textbf{H}=0$, equations (\ref{E1-1}) is reduced into $3$D incompressible Navier-Stokes equations, so the explicit blowup solutions (\ref{E1-6}) and (\ref{E1-6r1}) of Navier-Stokes equations should be a part of solutions for the corresponding MHD equations. 

\paragraph{One family of explicit blowup solutions.}

 When we set
\bel{X1-6}
\aligned
&v^r(t,r,z)={ar\over T^*-t},\\
&v^{\theta}(t,r,z)={k\over r},\\
&v^z(t,r,z)=-{2az\over T^*-t},
\endaligned
\ee
be a part of solutions for (\ref{X1-2})-(\ref{X1-2R4}), where constants $a,k\in\RR/\{0\}$.

Substituting (\ref{X1-6}) into equations (\ref{X1-2R0}) and (\ref{X1-2R2})-(\ref{X1-2R4}), we get
\bel{X1-2r0}
H^r\del_rH^{\theta}+H^z\del_zH^{\theta}+{1\over r}H^{\theta}H^r=0,
\ee
\bel{X1-2r2}
\del_tH^r+{ar\over T^*-t}\del_rH^r-{2az\over T^*-t}\del_zH^r=\nu(\triangle-{1\over r^2})H^r+{a\over T^*-t}H^r,
\ee
\bel{X1-2r3}
\del_tH^{\theta}+{ar\over T^*-t}\del_rH^{\theta}-{2az\over T^*-t}\del_zH^{\theta}+{2k\over r^2}H^r
=\nu(\triangle-{1\over r^2})H^{\theta}+{a\over T^*-t}H^{\theta},
\ee
\bel{X1-2r4}
\del_tH^z+{ar\over T^*-t}\del_rH^z-{2az\over T^*-t}\del_zH^z=\nu\triangle H^z-{2a\over T^*-t}H^z.
\ee

We observe the the structure of incompressibility condition (\ref{X1-4}) on the magnetic field, and we find it is better to set
\bel{E2-1}
\aligned
&H^r(t,r,z)={\bar{a}r\over (T^*-t)^{\alpha}},\\
&H^z(t,r,z)=-{2\bar{a}z\over (T^*-t)^{\alpha}},
\endaligned
\ee
where $\bar{a}\neq0,\alpha$ are two unknown constants.

It is easy to see $H^r(t,r,z)$ and $H^z(t,r,z)$ given in (\ref{E2-1}) satisfies the incompressibility condition (\ref{X1-4}) on the magnetic field.

Note that $(\triangle-{1\over r^2})r=0$.
Substituting the $H^r$ in (\ref{E2-1}) into (\ref{X1-2r2}), we get
$$
\alpha=0.
$$
which gives that 
\bel{E2-2}
\aligned
&H^r(t,r,z)=\bar{a}r,\\
&H^z(t,r,z)=-2\bar{a}z,
\endaligned
\ee

We now find $H^{\theta}(t,r,z)$. Assume that
\bel{E2-3}
H^{\theta}(t,r,z)={\bar{k}r^pz^q\over (T^*-t)^{\beta}},
\ee
where $\bar{k}\neq0$ and $p,q,\beta$ are unknown constants.

We substitute (\ref{E2-2})-(\ref{E2-3}) into (\ref{X1-2r0}) and (\ref{X1-2r3}), respectively, there are
\bel{E2-4}
p-2q+1=0,
\ee
and
\bel{E2-5}
{\beta\bar{k}r^pz^q\over (T^*-t)^{\beta+1}}+{ap\bar{k}r^pz^q\over (T^*-t)^{\beta+1}}-{2a\bar{k}qr^pz^q\over (T^*-t)^{\beta+1}}+{2k\bar{a}\over r}={a\bar{k}r^pz^q\over (T^*-t)^{\beta+1}}.
\ee

Equation (\ref{E2-5}) means that
$$
\beta=-1,\quad p=-1,\quad q=0,
$$
and 
$$
\bar{k}={2k\bar{a}\over 2a+1},\quad for \quad a\neq-{1\over2}.
$$
Thus we get
\bel{E2-6}
H^{\theta}={2k\bar{a}(T^*-t)\over (2a+1)r}.
\ee
Meanwhile, it is easy to check that (\ref{X1-2r4}) holds.

Hence by (\ref{E2-2}) and (\ref{E2-6}), we conclude 
$$
\aligned
&H^r(t,r,z)=\bar{a}r,\\
&H^{\theta}(t,r,z)={2k\bar{a}(T^*-t)\over (2a+1)r},\\
&H^z(t,r,z)=-2\bar{a}z,
\endaligned
$$
which combining with (\ref{X1-6}) gives a family of solutions of system (\ref{X1-2})-(\ref{X1-3}) with the incompressibility condition (\ref{X1-4}). Here constants $\bar{a},k\in\RR/\{0\}$ and $a\in\RR/\{-{1\over 2},0\}$.

Furthermore, we compute the pressure $P$. We substitute (\ref{X1-5}) into (\ref{X1-2}) and (\ref{X1-2R1}), there are
$$
\del_r\bar{P}=\Big(\bar{a}^2-{a(1+a)\over (T^*-t)^2}\Big)r-{k^2\Big(4\bar{a}^2(T^*-t)^2-(2a+1)^2\Big)\over r^3(2a+1)^2},
$$
and 
$$
\del_z\bar{P}=2z\Big(2\bar{a}^2+{a(1-2a)\over (T^*-t)^2}\Big).
$$
Note that 
$$
|\textbf{H}|^2=\bar{a}^2r^2+{4k^2\bar{a}^2(T^*-t)^2\over (2a+1)^2r^2}+4\bar{a}^2z^2.
$$
Thus by (\ref{X1-3}), direct computations give the pressure
$$
P(t,r,z)=-{1\over 2}\Big({a(a+1)r^2\over (T^*-t)^2}+{k^2\over r^2}\Big)+{az^2(1-2a)\over (T^*-t)^2}.
$$

\paragraph{Another family of explicit blowup solutions.}

When we set
\bel{X1-6r1}
\aligned
&v^r(t,r,z)={ar\over T^*-t},\\
&v^{\theta}(t,r,z)=kr(T^*-t)^{2a},\\
&v^z(t,r,z)=-{2az\over T^*-t},
\endaligned
\ee
be a part of solutions for (\ref{X1-2})-(\ref{X1-2R4}), where constants $a,k\in\RR/\{0\}$.

Substituting (\ref{X1-6r1}) into equations (\ref{X1-2R0}) and (\ref{X1-2R2})-(\ref{X1-2R4}), we get
\bel{X1-2rr0}
H^r\del_rH^{\theta}+H^z\del_zH^{\theta}+{1\over r}H^{\theta}H^r=0,
\ee
\bel{X1-2rr2}
\del_tH^r+{ar\over T^*-t}\del_rH^r-{2az\over T^*-t}\del_zH^r=\nu(\triangle-{1\over r^2})H^r+{a\over T^*-t}H^r,
\ee
\bel{X1-2rr3}
\del_tH^{\theta}+{ar\over T^*-t}\del_rH^{\theta}-{2az\over T^*-t}\del_zH^{\theta}+2k(T^*-t)^{2a}H^r
=\nu(\triangle-{1\over r^2})H^{\theta}+{a\over T^*-t}H^{\theta},
\ee
\bel{X1-2rr4}
\del_tH^z+{ar\over T^*-t}\del_rH^z-{2az\over T^*-t}\del_zH^z=\nu\triangle H^z-{2a\over T^*-t}H^z,
\ee
which compare with system (\ref{X1-2r0})-(\ref{X1-2r4}), we find there is only one equations different. This causes we have a new $H^{\theta}(t,r,z)$.
It is easy to check that $H^r(t,r,z)$ and $H^z(t,r,z)$ given in (\ref{E2-2})-(\ref{E2-3}) are the solution of (\ref{X1-2rr2}) and (\ref{X1-2rr4}) with
\bel{E2-4r0}
p-2q+1=0.
\ee
 
We substitute (\ref{E2-2})-(\ref{E2-3}) into (\ref{X1-2rr3}), there is
\bel{E2-5r0}
{\beta\bar{k}r^pz^q\over (T^*-t)^{\beta+1}}+{ap\bar{k}r^pz^q\over (T^*-t)^{\beta+1}}-{2a\bar{k}qr^pz^q\over (T^*-t)^{\beta+1}}+2k\bar{a}r(T^*-t)^{2a}={a\bar{k}r^pz^q\over (T^*-t)^{\beta+1}},
\ee
which gives that
$$
\alpha=-2a-1,\quad p=1,\quad q=1,
$$
and
$$
\bar{k}={2\bar{a}k\over 4a+1}.
$$
Thus we get
\bel{E2-6r0}
H^{\theta}(t,r,z)={2\bar{a}krz(T^*-t)^{2a+1}\over 4a+1}.
\ee

In conclusion, by (\ref{E2-2}) and (\ref{E2-6r0}), we obtain
$$
\aligned
&H^r(t,r,z)=\bar{a}r,\\
&H^{\theta}(t,r,z)={2\bar{a}krz(T^*-t)^{2a+1}\over 4a+1},\\
&H^z(t,r,z)=-2\bar{a}z,
\endaligned
$$
which combining with (\ref{X1-6}) gives another family of solutions of system (\ref{X1-2})-(\ref{X1-3}) with the incompressibility condition (\ref{X1-4}). Here constants $\bar{a},k\in\RR/\{0\}$ and $a\in\RR/\{-{1\over 4},0\}$.

Furthermore, we compute the pressure $P$. We substitute (\ref{X1-6}) into (\ref{X1-2}) and (\ref{X1-2R1}), there are
$$
\del_r\bar{P}=\Big(\bar{a}^2+k^2(T^*-t)^{4a}-{a(1+a)\over (T^*-t)^2}-{4\bar{a}^2k^2z^2(T^*-t)^{2(2a+1)}\over (4a+1)^2}\Big)r,
$$
and 
$$
\del_z\bar{P}=2z\Big(2\bar{a}^2+{a(1-2a)\over (T^*-t)^2}\Big).
$$
Note that 
$$
|\textbf{H}|^2=\bar{a}^2r^2+{4k^2\bar{a}^2r^2z^2(T^*-t)^{2(2a+1)}\over (4a+1)^2}+4\bar{a}^2z^2.
$$
Thus by (\ref{X1-3}), direct computations give the pressure
$$
\aligned
P(t,r,z)&={r^2\over 2}\Big(k^2(T^*-t)^{4a}-{a(a+1)\over (T^*-t)^2}-{8\bar{a}^2k^2z^2(T^*-t)^{2(2a+1)}\over (4a+1)^2}\Big)\\
&\quad+z^2\Big({a(1-2a)\over (T^*-t)^2}-{2\bar{a}^2k^2r^2(T^*-t)^{2(2a+1)}\over (4a+1)^2}\Big).
\endaligned
$$

\paragraph{Proof of Proposition 2.1.}
Since
$$
\aligned
&\textbf{v}(t,x)=v^r(t,r,z)\textbf{e}_r+v^{\theta}(t,r,z)\textbf{e}_{\theta}+v^z(t,r,z)\textbf{e}_z,\\
&\textbf{H}(t,x)=H^r(t,r,z)\textbf{e}_r+H^{\theta}(t,r,z)\textbf{e}_{\theta}+H^z(t,r,z)\textbf{e}_z,
\endaligned
$$
we can obtain a family of explicit blowup axisymmetric solutions for $3$D incompressible MHD equations by noticing that $\textbf{e}_r,\textbf{e}_{\theta},\textbf{e}_z$ are defined in (\ref{E1-0}), $r=\sqrt{x_1^2+x_2^2}$ and $z=x_3$, and
$$
\aligned
&v^r(t,r,z)={ar\over T^*-t},\\
&v^{\theta}(t,r,z)={k\over r},\\
&v^z(t,r,z)=-{2az\over T^*-t},\\
&H^r(t,r,z)=\bar{a}r,\\
&H^{\theta}(t,r,z)={2k\bar{a}(T^*-t)\over (2a+1)r},\\
&H^z(t,r,z)=-2\bar{a}z,
\endaligned
$$
where constants $\bar{a},k\in\RR/\{0\}$ and $a\in\RR/\{-{1\over2},0\}$, or 
$$
\aligned
&v^r(t,r,z)={ar\over T^*-t},\\
&v^{\theta}(t,r,z)=kr(T^*-t)^{2a},\\
&v^z(t,r,z)=-{2az\over T^*-t},\\
&H^r(t,r,z)=\bar{a}r,\\
&H^{\theta}(t,r,z)={2\bar{a}krz(T^*-t)^{2a+1}\over 4a+1},\\
&H^z(t,r,z)=-2\bar{a}z,
\endaligned
$$
where constants $\bar{a},k\in\RR/\{0\}$ and $a\in\RR/\{-{1\over 4},0\}$.

Futhermore the vorticity vector $\omega$ is 
$$
\omega(t,x)=\omega^r(t,r,z)\textbf{e}_r+\omega^{\theta}(t,r,z)\textbf{e}_{\theta}+\omega^z(t,r,z)\textbf{e}_z,
$$
where
$$
\aligned
&\omega^r(t,r,z)=-\del_zv^{\theta}=0,\\
&\omega^{\theta}(t,r,z)=\del_zv^r-\del_rv^z=0,\\
&\omega^z(t,r,z)={1\over r}\del_r(rv^{\theta})=0,
\endaligned
$$
or
$$
\aligned
&\omega^r(t,r,z)=-\del_zv^{\theta}=0,\\
&\omega^{\theta}(t,r,z)=\del_zv^r-\del_rv^z=0,\\
&\omega^z(t,r,z)={1\over r}\del_r(rv^{\theta})=2k(T^*-t)^{2a}.
\endaligned
$$

\paragraph{Proof of Theorem 1.1.}
By directly computations, we can obtain two family of explicit blowup solutions from (\ref{E1-6}) and (\ref{E1-6r1}).
Moreover, the vorticity vector 
$$
\omega(t,x)=\nabla\times\textbf{v}=0,
$$
or
$$
\omega(t,x)=\nabla\times\textbf{v}=2k(T^*-t)^{2a}.
$$




\

\

\

\textbf{Acknowledgments.} 

The author expresses his sincerely thanks to the BICMR of Peking University and Profes- sor Gang Tian for constant support and encouragement.
The author also expresses his sincerely thanks to Prof. V. Sver\'{a}k for informed the paper \cite{Con,K,K0} and his suggestions. 
The author is supported by NSFC No 11771359.


\begin{thebibliography}{xx}


\bibitem{A}
H. Alfv\'{e}n, Existence of electromagnetic-hydrodynamics waves. Nature. 150 (1942) 405-406.


\bibitem{CKN}
L. Caffarelli, R.V. Kohn, L. Nirenberg, Partial regularity of suitable weak solutions of the Navier-Stokes equations. Commun. Pure Appl. Math. XXXV (1982) 771-831.



\bibitem{CKS}
R.E. Caflisch, I. Klapper, G. Steele, Remarks on singularities, dimension and energy dissipation for ideal hydrodynamics and MHD. Comm. Math. Phys. 184 (1997) 443-455.



\bibitem{CL}
Y. Cai, Z. Lei, Global well-posedness of the incompressible Magnetohydrodynamics. Arch. Rational Mech. Anal. 228 (2018) 969-993.



\bibitem{Con}
P. Constantin, The Euler equations and non-local conservative Riccati equations. Internat. Math Res. Notices. 9 (2000) 455-465.





\bibitem{CM}
D. C\'{o}rdoba, C. Marliani, Evolution of current sheets and regularity of ideal in compressible magnetic fluids in 2D. Comm. Pure Appl. Math. 53 (2000) 512-524.


\bibitem{CP}
M. Cannone, F. Planchon,  Self-similar solutions for Navier-Stokes equations in $\RR^3$. Commun. Partial Differ. Equ. 21 (1996) 179-193.


\bibitem{Ch}
D. Chae, Nonexistence of self-similar singularities in the ideal Magnetohydrodynamics. Arch. Rational. Mech. Anal. 194 (2009) 1011-1027.



\bibitem{F}
C.L. Fefferman,  Existence and smoothness of the Navier-Stokes equations. Millenn. Prize Probl. (2006) 57-67.



\bibitem{GO}
J.D. Gibbon, K. Ohkitani, Evidence for singularity formation in a class of stretched solutions of the equations for ideal MHD, Tubes, sheets and singularities in fluid
dynamics (Zakopane, 2001). Fluid Mech. Appl. 71, 295-304 (2002)



\bibitem{HXY}
L.L. He, L. Xu, P. Yu, On Global dynamics of three dimensional Magnetohydro-dynamics: Nonlinear Stability of Alfv\'{e}n Waves (2016). arXiv:1603.08205.



\bibitem{HL}
T.Y. Hou, C. Li, Dynamic stability of the three-dimensional axisymmetric Navier-Stokes equations with
swirl. Commun. Pure Appl. Math. 61 (2008) 661-697.


\bibitem{HJL}
T. Y. Hou, T.L. Jin, P.F. Liu,
Potential singularity for a family of models of the axisymmetric incompressible flow. J. Nonlinear Sci. DOI 10.1007/s00332-017-9370-9.



\bibitem{JS}
H. Jia, V. Sver\'{a}k, Local in space estimates near initial time for weak solutions of the Navier-Stokes equations and forward self-similar solutions.
Invent. Math. 196 (2014) 233-265.

\bibitem{K}
L.V. Kapitanskii, Group analysis of the Navier-Stokes and Euler equations in the presence of rotation symmetry and new exact solutions to these equations. 
Dokl. Akad. Nauk SSSR 243 (1978), no. 4, 901-904.

\bibitem{K0}
L.V. Kapitanskii, Group analysis of the Navier-Stokes equations in the presence of rotational symmetry and some new exact solutions, Zap. Nauchn. Sem. LOMI 84 (1979) 89-107.





\bibitem{Lan}
 L.D. Landau, A new exact solution of the Navier-Stokes equations. C. R. (Doklady) Acad. Sci. URSS (N.S.) 43 (1944) 286-288.

\bibitem{Lan0}
L.D. Landau, E.M. Lifshitz, Fluid mechanics. Translated from the Russian by J.B. Sykes and W. H. Reid. Course of Theoretical Physics, Vol. 6. Pergamon Press, London-Paris-Frankfurt; Addison-Wesley Publishing Co., Inc., Reading, 1959



\bibitem{L}
J. Leray, Sur le mouvement d'un liquide visqueux emplissant l'espace. Acta Math. 63 (1934) 193-248.



\bibitem{Lei}
Z. Lei, On axially symmetric incompressible magnetohydrodynamics in three dimensions. J. Differential Equations. 259 (2015) 3202-3215.




\bibitem{LW}
J.G. Liu, W.C. Wang, Convergence analysis of the energy and helicity preserving scheme for axisymmetric
flows. SIAM J. Numer. Anal. 44 (2006) 2456-2480.



\bibitem{LH}
G. Luo, T.Y. Hou, Toward the finite-time blowup of the 3D incompressible Euler equations: a numerical
investigation. SIAM Multiscale Model. Simul. 12 (2014) 1722-1776.




\bibitem{NRS}
J. Necas, M. Ruzicka, V. Sver\'{a}k, On Leray's self-similar solutions of the Navier-Stokes equations. Acta Math. 176 (1996) 283-294.


\bibitem{Sv}
V. Sver\'{a}k, On Landau's solutions of the Navier-Stokes equations. J. Math. Sci. 179 (2011) 208-228.





\bibitem{Yan1}
W.P. Yan, Explicit self-similiar singularity of Born-Infeld equation, space-like surfaces with vanishing mean curvature equation and membrane equation.
arXiv: 1712.05159v3.



\bibitem{Yan2}
W.P. Yan, Two family of explicit blowup solutions for 3D incompressible Navier-Stokes equations. Preprint.




\bibitem{Yan0}
W.P. Yan, Finite time blowup explicit solution for 3D incompressible Euler equations. Preprint.












\end{thebibliography}
\end{document}